\documentclass[12pt, a4paper]{report}
\usepackage{amsfonts, amssymb, amsthm, latexsym, nccmath}  
\usepackage[cp1251]{inputenc} 
\usepackage[T2A]{fontenc}
\usepackage[english,russian,ukrainian]{babel}  
\usepackage{cite}
\usepackage{color, multirow} 
\usepackage{srcltx}
\usepackage{graphicx}
\usepackage{multicol}
\usepackage[centerlast]{caption2}
\usepackage{array,longtable}

\usepackage{wrapfig}
\graphicspath{{noiseimages/}}

\newtheorem{teorema}{Теорема}
\newtheorem{lema}{Лема}
\newtheorem{nas}{Наслідок}

\begin{document}

\noindent
{\bf С.О.~Чайченко}

\bigskip\noindent
{\bf S.O.~Chaichenko}

\bigskip\noindent
Україна, 84116, Донецька обл., м. Слов'янськ, вул. Генерала Батюка, 19

\bigskip\noindent
e-mail: s.chaichenko@gmail.com

\bigskip

\noindent
Державний вищий навчальний заклад "Донбаський державний педагогічний університет"

\bigskip

\bigskip\noindent
{\bf НАБЛИЖЕННЯ ЯДЕР БЕРГМАНА РАЦІОНАЛЬНИМИ \\ ФУНКЦІЯМИ З ФІКСОВАНИМИ ПОЛЮСАМИ}

\bigskip

\bigskip\noindent
{\bf APPROXIMATION OF THE BERGMAN KERNELS BY \\ RATIONAL FUNCTIONS WITH FIXED POLES}

\bigskip\bigskip

{\noindent \small Розв'язано задачу про найкращі раціональні наближення  ядер Бергмана  на одиничному колі комплексної площини у квадратичній і  рівномірній метриках.}

\bigskip
\bigskip

{\noindent \small Решено задачу о наилучших рациональных приближениях  ядер Бергмана  на единичной окруж\-ности комплексной плоскости в квадратичной и  равномерной метриках.}

\bigskip
\bigskip

{\noindent \small We solved the problem  of the best rational approximation of the Bergman kernels on the unit circle of the complex plane in the quadratic and uniform metrics.}

\bigskip
\bigskip

\noindent УДК 517.5

\newpage

У серії праць М.М.~Джрбашяна \cite{Djrbashan_1956, Djrbashan_1966-1, Djrbashan_1966-2, Dzhrbashyan-1967} було розвинено метод, що дозволяв отримувати розв'язки екстремальних задач про найкращі раціональні наближен\-ня ядра Коші на одиничному колі комплексної площини як у квадратичній, так і в рівномірній метриках. Цей метод, зокрема, спирався на зображення ядра Коші  за допомогою відрізку ряду Фур'є по ортонормованій на одиничному колі  системі раціональних функцій, що визначається фіксованою послідовністю полюсів, які лежать зовні  одиничного круга (система  Такенаки-Мальмквіста). У даній роботі, використовуючи зазначений метод і деякі результати роботи \cite{Dzhrbashyan-1967} знайдено розв'язок аналогічних задач для наближення вагових ядер Бергмана.

Наведемо означення і факти, які будемо використовувати в цій роботі.
Нехай ${\bf a}:=\{a_k\}_{k=0}^{\infty}$ --- послідовність точок в одиничному
крузі $\mathbb D:=\{z\in \mathbb{C}: |z|<1\}$ комплексної площини $\mathbb{C}$,
серед яких можуть бути точки скінченної і
навіть нескінченної кратності. Системою функцій
Такенаки-Мальмквіста, породженою послідовністю ${\bf a}$
називається  система
$\varphi:=\{\varphi_k\}_{k=0}^\infty$ функцій $\varphi_k$ вигляду \cite[\S 10.7]{Walsh_Interpol_1961}:
\begin{equation}\label{TM-system}
    \varphi_0(z):=\frac{\sqrt{1-|a_0|^2}}{1-\overline{a_0}z},~
    \varphi_k(z):=\frac{\sqrt{1-|a_k|^2}}{1-\overline{a_k}z}
    \prod_{j=0}^{k-1}\frac{-|a_j|}{a_j}\frac{z-a_j}{1-\overline{a_j}z},~
    k=1,2,\ldots,
\end{equation}
де при $a_j=0$ покладається $|a_j|/a_j=-1.$

Відомо \cite[\S 10.7]{Walsh_Interpol_1961}, що система Такенаки-Мальмквіста є
ортонормованою системою на одиничному колі $\mathbb T:=\{z\in \mathbb{C}: |z|=1\}$, тобто
\[
    \langle\varphi_k,\varphi_l\rangle:=\int_\mathbb{T} \varphi_k(z)\overline{\varphi_l(z)}~d\sigma(z) =\delta_{kl},\quad k,l=0,1,\ldots,
\]
де $\sigma$ -- нормована міра Лебега на $\mathbb T,$ $\delta_{kl}$ -- символ Кронекера.

Кожному елементу системи $\varphi$ поставимо у відповідність добуток
Бляшке $B_k$ степеня $k,$ тобто функцію вигляду
\[
    B_0(z):=1,\quad
    B_k(z):=\tau\prod_{j=0}^{k-1}\frac{z-a_j}{1-\overline{a_j}z},~n=1,2,\ldots,
\]
де $a_j\in\mathbb D,j=\overline{0,k-1},$ і $|\tau|=1.$

\begin{lema}\label{L-Hristoffel-equal}
Для довільних значень $z,\zeta \in \mathbb{D}$ і будь-якого $n \in \mathbb{N}$ справджується тотожність
\begin{equation}\label{Hristoffel-equal}
   \frac{1}{1-\overline z \zeta}=\sum_{k=0}^{n-1}
    \overline{\varphi_k(z)}\varphi_k(\zeta)+
    \frac{\overline{B_n(z)}{B_n(\zeta)}}{1-\overline z \zeta}.
\end{equation}
\end{lema}

Доведення леми \ref{L-Hristoffel-equal} можна знайти в \cite{Djrbashan_1956} (див. також роботу \cite{Savchuk_Best_linear_methods}).

Аналітичну у крузі $\mathbb{D}$ функцію $f$ відносять до простору Гарді $H_2(\mathbb{D}),$ якщо
\[\label{deff H2(D)}
    \|f\|_{ H_2}:= \sup_{0<\varrho<1}\left(\int_{\mathbb
    T}|f(\varrho z)|^2 d\sigma(z)\right)^{1/2}<\infty.
\]

Добре відомо, що функції з простору Гарді $H_2(\mathbb{D})$ мають на колі
$\mathbb T$ граничні значення по недотичних шляхах і
\begin{equation}\label{Cauchy-form}
    f(z)=\int_{\mathbb T}\frac{f(t)d\sigma(t)}{1-z\bar{t} },
    \quad\forall~z\in\mathbb D.
\end{equation}

З леми \ref{L-Hristoffel-equal} і формули (\ref{Cauchy-form}) випливає, що для будь-якої
функції $f \in H_2(\mathbb{D})$ і довільного $n \in \mathbb{N}$ справджується зображення \cite{Dzhrbashyan-1967}
\begin{equation}\label{represent}
    f(z)=\sum_{m=0}^{n-1} c_m(f) \varphi_m(z)+
    B_n(z)\int_\mathbb{T} \frac{\overline{B_n(t)} f(t)}{1-z\bar t }~d\sigma(t),
    \quad z \in \mathbb{D},
\end{equation}
де
$$
    c_m(f):=\int_\mathbb{T} f(t)\overline{\varphi_m(t)}~d\sigma(t),
    \quad  ~ m=0,1,\dots.
$$

У теорії просторів Бергмана (Bergman spaces) важливу роль відіграє функ\-ція вигляду (див., наприклад, \cite[p.~6]{Hedenmalm-TheorBergmanSp})
\begin{equation}\label{alpha-Bergman-kernel}
    \mathcal{K_\alpha}(z;w)= \frac{1}{(1-z\overline{w})^{2+\alpha}},
    \quad \alpha =0,1,\ldots, \quad z, w \in \mathbb{D},
\end{equation}
яку називають (ваговим) ядром Бергмана для круга $\mathbb{D}.$

Позначимо
\begin{equation}\label{Sn(K-alpha)}
    {S}_{n+1}(\mathcal{K_\alpha})(z;w):=\sum_{m=0}^{n} c_m (\mathcal{K_\alpha})\varphi_m (z)
\end{equation}
--- частинну суму ряду Фур'є функції $\mathcal{K_\alpha}(z;w)$ за  системою (\ref{TM-system}).

Справджується така лема.

\begin{lema}\label{L-rest-alphaBergKern}
Нехай $a_n=a_{n-1}=\ldots=a_{n-\alpha}\equiv w\in\mathbb{D}.$
Тоді для довільного $n=\alpha,\alpha+1,\ldots,$ і  $z\in\mathbb{D}$
\begin{equation}\label{rest-alphaBrgKern}
    \mathcal{K_\alpha}(z;w)-{S}_{n+1}(\mathcal{K_\alpha})(z;w)=
$$
$$
    =(-1)^\alpha
    \Big(\frac{\overline{ w}}{1-|w|^2}\Big)^{\alpha+1}
    \Bigg(\frac{w-z}{1-\overline{w}z}\Bigg)^{\alpha+1}\frac{B_{n-\alpha}(z) \overline{B_{n-\alpha}(w)}}{(\overline{w}z-1)}.
\end{equation}
\end{lema}

{\bf Доведення. } Оскільки при $w\in\mathbb{D}$ буде $\mathcal{K_\alpha}\in H_2(\mathbb{D}),$ то використовуючи зображення (\ref{represent}), отримуємо
\begin{equation}\label{represent-alphaBergKern}
    \mathcal{K_\alpha}(z;w)=\frac{1}{(1-z\overline{w})^{2+\alpha}}=
    \sum_{m=0}^{n} c_m(\mathcal{K_\alpha}) \varphi_m(z)+B_{n+1}(z) \mathcal{J}_{n,\alpha}(z;w),
\end{equation}
де
\begin{equation}\label{def-Jn-alpha}
    \mathcal{J}_{n,\alpha}(z;w):=\int_{\mathbb{T}}\frac{\overline{B_{n+1}(t)} \mathcal{K_\alpha}(t;w)}{1-z\bar{t}}d\sigma(t).
\end{equation}

Покладемо
$$
      \pi_n(z):=\prod_{m=0}^{n-\alpha-1}(a_m-z), \quad
      \tau_n(z):=\prod_{m=0}^{n-\alpha-1}(1-\bar{a}_m z), \quad n=\alpha, \alpha+1,\ldots~.
$$

Тоді
\begin{equation}\label{repres-Blashke-prod-alpha}
    B_{n+1}(z)=B_{n-\alpha}(z) \Big(\frac{w-z}{1-\overline{w}z}\Big)^{\alpha+1} \left(\frac{|w|}{w}\right)^{\alpha+1}=
$$
$$
     =\frac{\pi_{n-\alpha}(z)}{\tau_{n-\alpha}(z)} \prod_{m=0}^{n-\alpha-1} \frac{|a_m|}{a_m}
     \Big(\frac{w-z}{1-\overline{w}z}\Big)^{\alpha+1} \left(\frac{|w|}{w}\right)^{\alpha+1}.
\end{equation}
Звідси, використовуючи визначення функції $\mathcal{K_\alpha}(z;w)$ і очевидну тотожність
$$
     \overline{\Bigg( \frac{\pi_n (t)}{ \tau_n (t)}\Bigg)}=
     \frac{\tau_n (t)}{\pi_n (t)},  \quad t \in \mathbb{T},
$$
будемо мати
\begin{equation}\label{BnKalpha}
   \overline{B_{n+1}(t)} \mathcal{K_\alpha}(t;w)=
   \prod_{m=0}^{n-\alpha-1}\frac{|a_m|}{\overline{a}_m} \left(\frac{|w|}{\overline{w}}\right)^{\alpha+1}
   \Bigg(\frac{1-\overline{w}t }{w-t}\Bigg)^{\alpha+1}
   \frac{\tau_{n-\alpha} (t)}{\pi_{n-\alpha} (t)}\frac{1}{(1-\overline{w}t)^{2+\alpha}}=
$$
$$
    =\prod_{m=0}^{n-\alpha-1}\frac{|a_m|}{\overline{a}_m} \left(\frac{|w|}{\overline{w}}\right)^{\alpha+1}
    \frac{\tau_{n-\alpha} (t)}{\pi_{n-\alpha} (t)}
    \frac{1}{(w-t)^{\alpha+1}(1-\overline{w}t)}, \quad t \in \mathbb{T}.
\end{equation}

З рівностей (\ref{def-Jn-alpha}) i (\ref{BnKalpha}), одержуємо
$$
    \mathcal{J}_{n,\alpha}(z;w)=
    \prod_{m=0}^{n-\alpha-1}\frac{|a_m|}{\overline{a}_m}
    \left(\frac{|w|}{\overline{w}}\right)^{\alpha+1}
    \int_{\mathbb{T}} \frac{\tau_{n-\alpha}(t)}{\pi_{n-\alpha}(t)} \frac{1}{(w-t)^{\alpha+1}(1-\overline{w}t)}\frac{d\sigma(t)}{1-z\bar{t}}=
$$

\begin{equation}\label{repres-Jn-alpha}
    =\frac{1}{2\pi i}\prod_{m=0}^{n-\alpha-1}\frac{|a_m|}{\overline{a}_m} \left(\frac{|w|}{\overline{w}}\right)^{\alpha+1}
    \int_{\mathbb{T}} \frac{\tau_{n-\alpha}(t)}{\pi_{n-\alpha}(t)}
    \frac{1}{(w-t)^{\alpha+1}(t-z)} \frac{dt}{1-\overline{w}t}.
\end{equation}
У співвідношенні (\ref{repres-Jn-alpha}) підінтегральна функція
$$
    g_{n,\alpha}(t,w,z):=\frac{\tau_{n-\alpha}(t)}{\pi_{n-\alpha}(t)}
    \frac{1}{(w-t)^{\alpha+1}(t-z)}\frac{1}{1-\overline{w}t},
    \quad w,z \in \mathbb{D}, \quad \alpha=0,1,\ldots,
$$
як функція змінної $t,$ зовні одиничного круга $\mathbb{D}$ має єдиний простий полюс $1/\overline{w}$, а при $|t|\to \infty$ має порядок прямування до нуля $\mathcal{O}(|t|^{-2})$.
Ураховуючи ці факти і змінивши напрям інтегрування, з рівності (\ref{repres-Jn-alpha}) за теоремою про лишки отримуємо
$$
    \mathcal{J}_{n,\alpha}(z;w)=-\prod_{m=0}^{n-\alpha-1}\frac{|a_m|}{\overline{a}_m} \left(\frac{|w|}{\overline{w}}\right)^{\alpha+1}
    \lim_{t \to 1/\overline{w}} (t-1/\overline{w}) g_{n,\alpha}(t,w,z)=
$$
$$
   = \prod_{m=0}^{n-\alpha-1}\frac{|a_m|}{\overline{a}_m} \left(\frac{|w|}{\overline{w}}\right)^{\alpha+1}
   \frac{ 1/\overline{w} }{(w-1/\overline{w})^{\alpha+1}(1/\overline{w}-z)}
    \frac{\tau_{n-\alpha} (1/\overline{w})}{\pi_{n-\alpha}(1/\overline{w})}=
$$
\begin{equation}\label{equal-J_n-alpha}
    =(-1)^\alpha \Big(\frac{|w|}{1-|w|^2}\Big)^{\alpha+1} \frac{1}{\overline{w}z-1}
    \overline{\left(\frac{\prod\limits_{m=0}^{n-\alpha-1}(a_m-w)}
    {\prod\limits_{m=0}^{n-\alpha-1}(1-\overline{a}_m w)} \right)}
    \prod_{m=0}^{n-\alpha-1} \frac{|a_m|}{\overline{a}_m}=
$$
$$
    =(-1)^\alpha \Big(\frac{|w|}{1-|w|^2}\Big)^{\alpha+1}
    \frac{\overline{B_{n-\alpha}(w)}}{\overline{w}z-1}.
\end{equation}

Оскільки $|w|^2/w=\overline{w},$ то з рівностей (\ref{represent-alphaBergKern}), (\ref{repres-Blashke-prod-alpha}) і (\ref{equal-J_n-alpha}) випливає формула (\ref{rest-alphaBrgKern}).

Лему доведено.

Для довільного фіксованого значення $\omega\in\mathbb{D}$ визначимо раціональну функ\-цію
\begin{equation}\label{rn-funct}
    r_{\alpha,n}(z;w)=(1-z\overline{w}) S_{n+1}(\mathcal{K_\alpha})(z;w).
\end{equation}

Справджується таке твердження.

\begin{lema}\label{lema2}
Нехай $a_n=a_{n-1}=\ldots=a_{n-\alpha}\equiv w\in\mathbb{D}.$ Тоді для довільного $n\in \mathbb{N}$ і $z\in\mathbb{D}$
\begin{equation}\label{rn-equal}
      r_{\alpha,n}(z;w)=\frac{(1-|w|^2)^{\alpha+1}+(-1)^\alpha(\overline{w}(w-z))^{\alpha+1}
      B_{n-\alpha}(z)\overline{B_{n-\alpha}(w)}}{(1-|w|^2)^{\alpha+1}(1-\overline{w}z)^{\alpha+1}}
\end{equation}
і функція $r_{\alpha,n}(z;w)$ задовольняє інтерполяційні умови
\begin{equation}\label{rn-iterp-condit}
      r_{\alpha,n}^{(s_m-1)}(a_m;w)=
      \frac{(\alpha +s_m-1)!~\overline{w}^{s_m-1}}{(1-\overline{w}a_m)^{\alpha+s_m}}, \quad
      0\leq m\leq n.
\end{equation}
       \end{lema}
{\bf Доведення.} Помноживши ліву і праву частину формули (\ref{rest-alphaBrgKern}) на $(1-\overline{w}z)$
і врахувавши співвідношення (\ref{rn-funct}), отримаємо рівність (\ref{rn-equal}).

Нехай $s_k \geq 1$ --- кратність появи числа $a_m$
в послідовності $\{a_j\}_{j=0}^m$, а $p_m(n)$ --- кратність появи числа $a_m$
в $\{a_j\}_{j=0}^n$. Зрозуміло, що $1\leq s_m \leq p_m (n)$, $0\leq m\leq n$ i
$s_n=p_n(n)$. Отже, якщо $0 \leq j\le s_k-1 \leq p_k (n),$  то функція
$$
    \psi_{nj}(t)=\frac{\pi_n (t)}{(t-a_m)^{j+1}}
$$
є поліномом. Диференціюючи $j$ разів по $z$  формулу
\begin{equation}\label{pi-n-int}
      \pi_n(z)=\frac{1}{2\pi i} \int_\mathbb{T}\frac{\pi_n (t)}{t-z}dt, \quad
      |z|\leq 1,
\end{equation}
знаходимо
$$
    \frac{d^j}{dz^j}[\pi_n(z)]=\frac{j!}{2\pi i}\int_\mathbb{T}
    \frac{\pi_n(t)}{(t-z)^{j+1}} dt,
$$
звідки при $z=a_m$ одержуємо рівність
\begin{equation}\label{deriv(pi-n)}
      [\pi_n(a_m)]^{(j)}=\frac{j!}{2\pi i} \int_\mathbb{T}\frac{\pi_n(t)}{(t-a_m)^{j+1}}dt=
      \frac{j!}{2\pi i}\int_\mathbb{T}\psi_{nj}(t) dt=0,
\end{equation}
яка виконується при всіх $0\leq j\leq s_m-1.$

Далі зауважимо, що
\begin{equation}\label{deriv-Coushy-kern}
     \frac{d^{s_m-1}}{dz^{s_m-1}}\Bigg(\frac{1}{(1-z\overline{w})^{\alpha+1}}\Bigg)=
     \frac{(\alpha+s_m-1)!~\overline{w}^{s_m-1}}{(1-z\overline{w})^{\alpha+s_m}},
     \quad 0\le m\le n.
\end{equation}

Нарешті, після $(s_m-1)$-кратного диференціювання по $z$ тотожності
$$
    \frac{1}{(1-z\overline{w})^{\alpha+1}}-r_{\alpha,n}(z;w)
     =(-1)^\alpha\Bigg(\frac{ |w|}{1-|w|^2}\Bigg)^{\alpha+1} 
     \overline{B_{n-\alpha}(w)} B_{n+1}(z)=
$$
$$
    =(-1)^\alpha \Bigg(\frac{ |w|}{1-|w|^2}\Bigg)^{\alpha+1} \overline{B_{n-\alpha}(w)}~\frac{\pi_{n+1}(z)}{\tau_{n+1}(z)} 
    \prod_{m=0}^{n-1} \frac{|a_m|}{a_m}, 
$$
у якій $a_n=\ldots=a_{n-\alpha}\equiv w\in\mathbb{D},~z\in\mathbb{D},$
на підставі співвідношень (\ref{deriv(pi-n)}) і (\ref{deriv-Coushy-kern}), отримуємо (\ref{rn-iterp-condit}),  оскільки
$$
      \frac{(\alpha+s_m-1)!~\overline{w}^{s_m-1}}{(1-a_m\overline{w})^{1+\alpha}}-
      r_{\alpha,n}^{(s_m-1)}(a_m;w)=
$$
$$
    =(-1)^\alpha\Bigg(\frac{ |w|}{1-|w|^2}\Bigg)^{\alpha+1} \overline{B_{n-\alpha}(w)} 
    \prod_{k=0}^{n-1} \frac{|a_k|}{a_k} \sum_{j=0}^{s_m-1} C_{s_m-1}^j \frac{\pi_{n+1}^{(j)}(a_m)}{\tau_{n+1}^{(s_m-j)}(a_m)}=0.
$$

Лему доведено.

Для даного $n$ $(1\leq n<\infty)$ позначимо через  $\mathcal{R}(n)$ множину раціональних функцій вигляду
$$
    R_n(x)=c_0+\sum_{m=1}^{n}c_m\varphi_{m-1}(x).
$$

Беручи до уваги визначення системи $\{\varphi_m(z)\}_{m=0}^\infty$, легко зрозуміти, що $\mathcal{R}(n)$ збігається з множиною раціональних функцій вигляду
$$
    b_0^{(n)}+\sum_{j=1}^{n}\frac{b_j^{(n)}}{(1-\bar{a}_jx)^{s_j}}.
$$

Основним результатом роботи є таке твердження.

\begin{teorema}\label{T-inf-mu}
На множині функцій $\mathcal{R}(n)$ мінімум функціоналу
\begin{equation}\label{mu-functional-alpha}
    \mu_\alpha ( R_n):=\int_\mathbb{T} \Big|\frac{1}{(1-x\overline{w})^{2+\alpha}}-\frac{R_n(x)}{(1-x\overline{w})}\Big|^2 d\sigma(x)
\end{equation}
реалізує функція
$$
    r_{\alpha,n}(x;w)=(1-x\overline{w})S_{n+1}(\mathcal{K_\alpha})(x;w)\in \mathcal{R}(n),
$$
де $S_{n+1}(\mathcal{K_\alpha})(x;w)$ --- частинна сумма порядку $n+1$ ряду Фур'є функції $\mathcal{K_\alpha}(x;w)$  по системі (\ref{TM-system}), $a_n=a_{n-1}=\ldots=a_{n-\alpha}\equiv w\in\mathbb{D},$  і при цьому виконується рівність
\begin{equation}\label{inf-mu-alpha}
    \inf_{R_n\in \mathcal{R}(n)} \mu_\alpha(R_n)=
    \mu_\alpha(r_{\alpha,n})=\frac{|w |^{2\alpha+2}}{(1-|w|^2)^{2\alpha+3}} |B_{n-\alpha}(w)|^2.
\end{equation}
\end{teorema}

{\bf Доведення.} З тотожності (\ref{rest-alphaBrgKern}) леми \ref{L-rest-alphaBergKern} випливає рівність
\begin{equation}\label{gener-equal}
    \int_\mathbb{T} \Big|\frac{1}{(1-x\overline{w})^{2+\alpha}}-S_{n+1}(\mathcal{K_\alpha})(x;w)\Big|^2 d\sigma(x)=
$$
$$
    =\Big(\frac{ |\overline{w} | }{(1-|w|^2)}\Big)^{2\alpha+2} |B_{n-\alpha} (w)|^2
    \int_\mathbb{T} \frac{|B_{n+1}(x)|^2}{|1-x\overline{w}|^2}d\sigma(x)=
$$
$$
   =\frac{|w |^{2\alpha+2}}{(1-|w|^2)^{2\alpha+3}} |B_{n-\alpha}(w)|^2.
\end{equation}

Зауважимо, що довільна раціональна функція $R_n(x)\in \mathcal{R}(n)$
дозволяє зображення вигляду
\begin{equation}\label{Rn-repres}
    R_n(x)=(1-x\overline{w})\sum_{m=0}^n c_m\varphi_m(x), \quad a_n=\ldots=a_{n-\alpha}\equiv w\in\mathbb{D},
\end{equation}
з певними коефіцієнтами $\{c_m\}_{k=0}^n$ і навпаки, для довільного набору коефіцієнтів
$\{c_m\}_{m=0}^n$ вираз з правої частини рівності (\ref{Rn-repres}) є деякою раціональною функцією
$R_n\in \mathcal{R}(n)$.

Внаслідок цього і беручи до уваги визначення функції $r_{\alpha,n}(x;w)$ (рівність (\ref{rn-funct})),
можна стверджувати, що для довільної функції  $R_n\in \mathcal{R}(n)$ функціонал
$\mu_\alpha(R_n)$  має зображення
\begin{equation}\label{mu-repres}
        \mu_\alpha(R_n)=\int_\mathbb{T}\Big|\frac{1}{(1-x\overline{w})^{2+\alpha}}-
        \sum_{m=0}^n c_m\varphi_m(x)\Big|^2 d\sigma (x),
        \end{equation}
з певними коефіцієнтами $c_m,~m=0,1,\ldots,n,$ і навпаки, для довільного набору сталих $c_m,~m=0,1,\ldots,n,$ вираз з правої частини співвідношення (\ref{mu-repres})  є значенням функціоналу
$\mu_\alpha(R_n)$  при деякому $R_n\in \mathcal{R}(n)$.

Ураховуючи те, що сума $S_{n+1}(\mathcal{K_\alpha})(x;w)$ є
$(n+1)$-м відрізком ряду Фур'є функції $\mathcal{K_\alpha}(x;w)$ по ортонормованій на одиничному колі $\mathbb{T}$ системі $\{\varphi_m\}_{m=0}^\infty,$ зі співвідношень (\ref{gener-equal}) і (\ref{mu-repres})   для довільної раціональної функції $R_n\in \mathcal{R}(n)$ отримуємо
\begin{equation}\label{ivers-inequal-mu}
    \mu_\alpha(R_n) \ge \mu_\alpha (r_{\alpha,n})=
$$
$$
    =\int_\mathbb{T} \Big| \frac{1}{(1-x\overline{w})^{2+\alpha}}-S_{n+1}(\mathcal{K_\alpha})(x;w)\Big|^2 d\sigma(x)=
    \frac{|w|^{2\alpha+2}}{(1-|w|^2)^{2\alpha+3}} |B_{n-\alpha}(w)|^2.
\end{equation}

При цьому,  знак рівності у (\ref{ivers-inequal-mu})  можливий лише у випадку, коли
$$
    R_n(x)=r_{\alpha,n}(x;w), \quad w\in\mathbb{D}.
$$

Співставлення співвідношень (\ref{gener-equal}) і (\ref{ivers-inequal-mu}) переконує нас у справедливості формули (\ref{inf-mu-alpha}).

Теорему доведено.

У випадку $\alpha=0$  ядра (\ref{alpha-Bergman-kernel}) мають вигляд
$$
    \mathcal{K}_0(z;w):=\mathcal{K}(z;w)= \frac{1}{(1-z\overline{w})^2}, \quad z, w \in \mathbb{D},
$$
і називаються (звичайними) ядрами Бергмана \cite[p.~6]{Hedenmalm-TheorBergmanSp}.

Покладаючи $\alpha=0,$ з теореми \ref{T-inf-mu} отримуємо такий наслідок.

\begin{nas}\label{Nas-inf-mu}
На множині функцій $\mathcal{R}(n)$ мінімум функціоналу
\begin{equation}\label{mu-functional}
    \mu (\mathcal{K}; R_n):=
    \int_\mathbb{T} \Big|\frac{1}{(1-x\overline{w})^2}-\frac{R_n(x)}{(1-x\overline{w})}\Big|^2 d\sigma(x)
\end{equation}
реалізує функція
$$
    r_n(x;w)=(1-x\overline{w})S_{n+1}(\mathcal{K})(x;w)\in \mathcal{R}(n),
$$
де $S_{n+1}(\mathcal{K})(x;w)$ --- частинна сумма порядку $n+1$ ряду Фур'є функції $\mathcal{K}(x; w)$  по системі (\ref{TM-system}), $a_n\equiv w\in\mathbb{D},$  і при цьому виконується рівність
\begin{equation}\label{inf-mu}
    \inf_{R_n\in \mathcal{R}(n)} \mu(\mathcal{K};R_n)=
    \mu(\mathcal{K};r_n)=\frac{|w B_n(w)|^2}{(1-|w|^2)^3}.
\end{equation}
\end{nas}

Використовуючи лему \ref{L-rest-alphaBergKern} i теорему \ref{T-inf-mu} отримуємо таке твердження.

\begin{teorema}\label{T-Coushy-kern}
Серед раціональних функцій $R_n\in \mathcal{R}(n)$ мінімум функціоналу
\begin{equation}\label{def-nu-alpha}
    \nu_\alpha (R_n):=\sup_{|x|=1} \Big| \frac{1}{(1-x\overline{w})^{1+\alpha}}-R_n(x) \Big|,
    \quad  R_n\in \mathcal{R}(n), \quad \alpha=0,1,\ldots,
\end{equation}
реалізує  функція
$$
    r_{\alpha,n}(x;w)=(1-x\overline{w})S_{n+1}(\mathcal{K_\alpha})(x;w)\in \mathcal{R}(n),
$$
де $S_{n+1}(\mathcal{K_\alpha})(x;w)$ --- частинна сумма порядку $n+1$ ряду Фур'є функції $\mathcal{K_\alpha}(x;w)$ по системі (\ref{TM-system}), $a_n=a_{n-1}=\ldots=a_{n-\alpha}\equiv w\in\mathbb{D},$  i при цьому
\begin{equation}\label{inf-nu-alpha}
    \inf_{R_n\in \mathcal{R}(n)} \nu_\alpha (R_n)=
    \nu_\alpha (r_{\alpha,n})=\Bigg(\frac{|w |}{1-|w|^2}\Bigg)^{1+\alpha} |B_{n-\alpha}(w)|.
\end{equation}
\end{teorema}

Зі співвідношення (\ref{def-nu-alpha}) випливає, що при $\alpha=0$ величина з правої частини рівності (\ref{inf-nu-alpha}) є найкращим рівномірним наближенням ядра Коші, тобто функції вигляду
$$
    \mathcal{C}(z;w):=\frac{1}{1-z\overline{w}}, \quad z,w \in \mathbb{D},
$$
на одиничному колі за допомогою раціональних функцій з множини $\mathcal{R}(n).$
Її точне значення вперше було отримане в \cite{Savchuk_Best_linear_methods}.

{\bf Доведення.}   Дійсно, зі співвідношень (\ref{mu-functional-alpha}), (\ref{inf-mu-alpha}) і (\ref{inf-nu-alpha}) одержуємо
$$
    \Big(\frac{|w |}{1-|w|^2}\Big)^{2\alpha+2} \frac{|B_{n-\alpha}(w)|^2}{1-|w|^2}
    \le \mu_\alpha (R_n)=
    \int_\mathbb{T} \Big|\frac{1}{(1-x\overline{w})^{2+\alpha}}-\frac{R_n(x)}{1-x\overline{w}} \Big|^2 d\sigma(x)=
$$
$$
    =\int_\mathbb{T} \Big |\frac{1}{(1-x\overline{w})^{1+\alpha}}-R_n(x) \Big|^2\frac{d\sigma(x)}{|1-x\overline{w}|^2} \le
    \nu^2_\alpha(R_n) \int_\mathbb{T} \frac{d\sigma(x)}{|1-x\overline{w}|^2}
    =\frac{\nu^2_\alpha(R_n)}{1-|w|^2},
$$
звідки випливає, що для довільного $R_n\in \mathcal{R}(n)$
$$
    \inf_{R_n\in w(a_k)} \nu_\alpha (R_n)\geq
    \Bigg(\frac{|w |}{1-|w|^2}\Bigg)^{1+\alpha} |B_{n-\alpha}(w)|.
$$

З іншого боку, використовуючи тотожність (\ref{rest-alphaBrgKern}), знаходимо
$$
    \nu_\alpha ( r_{\alpha,n})= \sup_{x\in \mathbb{T}}
    \Big|\frac{1}{(1-x\overline{w})^{1+\alpha}}-
    (1-x\overline{w}) S_{n+1}(\mathcal{K_\alpha})(x;w)\Big|=
$$
$$
     =\sup_{x\in \mathbb{T}}|1-x\overline{w}|
     \Big|\frac{1}{(1-x\overline{w})^{2+\alpha}}-S_{n+1}(\mathcal{K_\alpha})(x;w)\Big|=
$$
$$
     =\sup_{x\in\mathbb{T}}|1-x\overline{w}| \bigg|\bigg(\frac{\overline{w}}{1-|w|^2}\bigg)^{1+\alpha}
      \overline{B_{n-\alpha}(w)} \frac{B_{n+1}(x)}{\overline{w}x-1}\bigg|=
$$
$$
     =\Bigg(\frac{|w |}{1-|w|^2}\Bigg)^{1+\alpha} |B_{n-\alpha}(w)| \sup_{x\in\mathbb{T}} |B_{n+1}(x)|=
     \Bigg(\frac{|w |}{1-|w|^2}\Bigg)^{1+\alpha}|B_{n-\alpha}(w)|.
$$

Теорему доведено.

\bibliographystyle{plain}
\renewcommand{\refname}{References}

\end{document}